\documentclass[a4paper,11pt]{article}
\usepackage{times}
\usepackage{enumerate}
\def\titlerunning#1{\gdef\titrun{#1}}
\makeatletter
\def\author#1{\gdef\autrun{\def\and{\unskip, }#1}\gdef\@author{#1}}
\def\address#1{{\def\and{\\\hspace*{18pt}}\renewcommand{\thefootnote}{}%
\footnote {#1}}%
\markboth{\autrun}{\titrun}}
\makeatother

\usepackage{amsmath}
\usepackage{amsfonts}
\usepackage{amssymb}
\usepackage[english]{babel}
\usepackage{graphicx}
\usepackage{amsthm}
\usepackage{xcolor}

\usepackage{wasysym}

\usepackage{bbm}

\theoremstyle{definition}

\theoremstyle{remark}

\numberwithin{equation}{section}

\newcommand{\Z}{\mathbb{Z}}

\begin{document}
\baselineskip=17pt
\titlerunning{}
\title{The first example of a simple $2-(81,6,2)$ design}

\author{Anamari Nakic}

\date{}

\maketitle

\address{Faculty of Electrical Engineering and Computing,
University at Zagreb, Croatia, email: anamari.nakic@fer.hr}

\begin{abstract}
\noindent
We give the very first example of a simple $2-(81,6,2)$ design. Its points are the elements of the
elementary abelian group of order 81 and each block is the union of two parallel lines of the
4-dimensional geometry over the field of order 3. Hence it is also additive. 
\end{abstract}


\section{Introduction}	

A $t-(v,k,\lambda)$ design is a pair $(V,{\cal B})$ where  $V$ is a set of $v$ {\it points} and ${\cal B}$ is a collection of $k$-subsets ({\it blocks}) of $V$ with the property that every $t$-subset of $V$ is contained in exactly $\lambda$ blocks.
A design is said to be {\it simple} if it does not have repeated blocks, i.e., if the multiset of blocks
is actually a set.
We recall, in particular, that in a $2-(v,k,\lambda)$ design the number of blocks
containing any fixed point is ${\lambda (v-1)\over k-1}$ and that the number of blocks is ${\lambda v(v-1)\over k(k-1)}$.
Thus $\lambda(v-1)\equiv0$ (mod $k-1$) and $\lambda v(v-1)\equiv0$ (mod $k(k-1)$) are trivial necessary
conditions (called {\it divisibility conditions}) for the existence of a $2-(v,k,\lambda)$ design.

Two designs $(V,{\cal B})$ and $(V,'{\cal B}')$ are {\it isomorphic} if there
exists a bijection between $V$ and $V'$ turning ${\cal B}$ into ${\cal B}'$.
An automorphism group of a design $(V,{\cal B})$ is a group $A$ of permutations on $V$
leaving ${\cal B}$ invariant. It is convenient to have designs with a ``rich" automorphism group $A$
since they can be stored more efficiently in terms of memory space; we do not need to give the entire 
collection of blocks but only a system of representatives for the $A$-orbits on it.

For general background on the topic we refer to \cite{BJL}.

The authors of \cite{CFP} developed an interesting theory about {\it additive} designs. These are designs
$(V,{\cal B})$ for which, up to isomorphism, $V$ is a subset of a suitable additive group $G$ and the elements 
of any block sum up to zero. We propose to speak of a {\it strictly} additive design if $V$ {\it coincides} with the set of elements of $G$. 

Among the strictly additive 2-designs with $\lambda=1$ we obviously have the point-line $2-(q^n,q,1)$ designs associated with 
AG$(n,q)$, the $n$-dimensional affine geometry over the field of order $q$. As far as we are aware, no other example of a 
strictly additive $2-(v,k,1)$ design is known. In the attempt to find strictly additive 2-designs with a higher $\lambda$
it is natural to look for 2-$(q^n,mq,\lambda)$ designs whose points are those of
AG$(n,q)$ and whose blocks are union of $m$ pairwise disjoint lines. Indeed such a
design would be strictly additive automatically. It is quite evident that the set
of all possible unions of two parallel lines of an affine plane AG$(2,q)$ is a $2-(q^2,2q,2q-1)$ design.
For $q$ odd, this design has been characterized in \cite{CF} as the unique design with these parameters and the
property that the intersection of any two distinct blocks has size 0, 4 or $q$. 

In this note we explicitly give a simple $2-(81,6,2)$ design which is strictly additive.
Indeed its points are those of AG$(4,3)$ and each block is the union of two parallel lines.
As far as we are aware, this is the very first example of a simple design with these parameters.

\section{The design}

The values of $v$ for which a $2-(v,6,\lambda)$ design exists have been completely 
determined for $\lambda>1$ by Hanani \cite{H}. For $\lambda=1$ the existence is still uncertain 
for twenty-nine values of $v$ the third of which is $81$. Thus it is unknown whether a
$2-(81,6,1)$ design exists. According to the tables of 2-designs of small order by Mathon and Rosa \cite{MR}
there is only one known $2-(81,6,2)$ design. Looking at its explicit description in Examples 2.7 on page 237
of \cite{GM}, any expert reader will recognize that it has sixteen pairs of repeated blocks. Thus it cannot be 
isomorphic to the simple and strictly additive design with the same parameters that we explicitly construct below.

Let $G=\Z_3^4$ be the elementary abelian group of order 81.  Given two elements $x\in G\setminus\{0\}$ 
and $y\in G\setminus\{0,x,2x\}$, let
$B(x,y)$ be the union of the two parallel lines $\{0,x,2x\}$ and $\{y,x+y,2x+y\}$ of AG$(4,3)$. 
The $G$-stabilizer of $B(x,y)$ (under the natural action of $G$ on itself) is clearly given
by $\{0,x,2x\}$, hence its $G$-orbit has size $|G|/3=27$. 
Also, from the divisibility conditions we infer that a $2-(81,6,2)$ design has $432=27\cdot16$ blocks.
Thus it makes sense to look for a design with these parameters whose collection of blocks
is the union of the $G$-orbits of 16 suitable blocks of the form $B(x,y)$.
Such a 16-tuple of blocks has been found  with a computer and it is given below.

$$\{(0, 0, 0, 0), (0, 0, 0, 1), (0, 0, 0, 2), (0, 1, 0, 0), (0, 1, 0,  1), (0, 1, 0, 2)\}$$
$$\{(0, 0, 0, 0), (0, 0, 1, 1), (0, 0, 2, 2), (2, 1, 0, 0), (2, 1, 1, 1), (2, 1, 2, 2)\}$$
$$\{(0, 0, 0, 0), (0, 1, 1, 1), (0, 2, 2, 2), (0, 0, 1, 0), (0, 1, 2, 1), (0, 2, 0, 2)\}$$
$$\{(0, 0, 0, 0), (0, 1, 2, 0), (0, 2, 1, 0), (2, 0, 2, 1), (2, 1, 1, 1), (2, 2, 0, 1)\}$$
$$\{(0, 0, 0, 0), (1, 0, 0, 0), (2, 0, 0, 0), (0, 2, 2, 1), (1, 2, 2, 1), (2, 2, 2, 1)\}$$
$$\{(0, 0, 0, 0), (1, 0, 1, 0), (2, 0, 2, 0), (0, 1, 0, 0), (1, 1, 1, 0), (2, 1, 2, 0)\}$$
$$\{(0, 0, 0, 0), (1, 0, 1, 1), (2, 0, 2, 2), (0, 0, 2, 0), (1, 0, 0, 1), (2, 0, 1, 2)\}$$
$$\{(0, 0, 0, 0), (1, 0, 2, 0), (2, 0, 1, 0), (0, 2, 1, 1), (1, 2, 0, 1), (2, 2, 2, 1)\}$$
$$\{(0, 0, 0, 0), (1, 0, 2, 2), (2, 0, 1, 1), (0, 1, 2, 1), (1, 1, 1,  0), (2, 1, 0, 2)\}$$
$$\{(0, 0, 0, 0), (1, 1, 0, 0), (2, 2, 0, 0), (0, 2, 0, 1), (1, 0, 0, 1), (2, 1, 0, 1)\}$$
$$\{(0, 0, 0, 0), (1, 1, 0, 1), (2, 2, 0, 2), (0, 2, 2, 0), (1, 0, 2, 1), (2, 1, 2, 2)\}$$
$$\{(0, 0, 0, 0), (1, 1, 2, 0), (2, 2, 1, 0), (0, 0, 2, 1), (1, 1, 1, 1), (2, 2, 0, 1)\}$$
$$\{(0, 0, 0, 0), (1, 1, 2, 1), (2, 2, 1, 2), (0, 2, 1, 1), (1, 0, 0, 2), (2, 1, 2, 0)\}$$
$$\{(0, 0, 0, 0), (1, 1, 2, 2), (2, 2, 1, 1), (0, 2, 2, 0), (1, 0, 1, 2), (2, 1, 0, 1)\}$$
$$\{(0, 0, 0, 0), (1, 2, 1, 2), (2, 1, 2, 1), (0, 0, 2, 1), (1, 2, 0, 0), (2, 1, 1, 2)\}$$
$$\{(0, 0, 0, 0), (1, 2, 2, 0), (2, 1, 1, 0), (0, 2, 2, 1), (1, 1, 1, 1), (2, 0, 0, 1)\}$$

Here is a short program in GAP \cite{GAP} checking that the union of the $G$-orbits of the above sixteen  
6-subsets of $G$ actually is the collection of blocks of the desired $2-(81, 6, 2)$ design.

\texttt{\# All points of AG(4,3)}

\texttt{pts := Tuples( [0..2],4 );}

\texttt{\# blkOrbRep denotes block orbit representatives listed above}

\texttt{\# all blocks of the design}

\texttt{blks:=Union(List (blkOrbRep , b-> List(pts, p -> }

\texttt{AsSet(List([1..Size(b)], i-> (b[i] + p) mod 3)) ) ));;}

\texttt{\# check that it is a 2-design}

\texttt{Collected(List(Combinations(pts, 2), p -> }

\texttt{Number(blks, b-> (p[1] in b) and (p[2] in b) ) ));}

\medskip
It is evident that any block of the obtained design is a union of two parallel lines.
Hence we conclude that this design is strictly additive. It is also easy to check 
that our design is simple. Thus, considering the comments that we made at the beginning of this section, 
we have the following new result.

\begin{quote}
The number of pairwise non-isomorphic $2-(81,6,2)$ designs is at least equal to $2$.
\end{quote}

Some infinite classes of strictly additive 2-designs will be given in a future paper 
still in preparation \cite{BN}.

\section*{Acknowledgements}
The author is supported by the Croatian Science Foundation under the project 9752.

\end{document}